\documentclass[letterpaper, 10 pt, conference]{ieeeconf}
\usepackage{graphicx}

\usepackage{amsmath}
\usepackage{amsfonts}
\newtheorem{defin}{Definition}
\newtheorem{theorem}{Theorem}

\newcommand{\pder}[2][]{\frac{\partial#1}{\partial#2}}

\newcommand{\RR}{\mathbb{R}}

\begin{document}

\title{\LARGE \bf Output-Feedback Control of Nonlinear Systems using Control Contraction Metrics and Convex Optimization}

\IEEEoverridecommandlockouts
\author{
Ian R. Manchester$^1$ \ \ Jean-Jacques E. Slotine$^2$ \\ \ \\
1: ACFR, School of Aerospace, Mechanical and Mechatronic Engineering, University of Sydney, Australia\\
2: Nonlinear Systems Laboratory, Massachusetts Institute of Technology, USA\\
ian.manchester@sydney.edu.au \ jjs@mit.edu}
\maketitle

\begin{abstract}
Control contraction metrics (CCMs) are a new approach to nonlinear control design based on contraction theory. The resulting design problems are expressed as pointwise linear matrix inequalities and are and well-suited to solution via convex optimization. In this paper, we extend the theory on CCMs by showing that a pair of ``dual'' observer and controller problems can be solved using pointwise linear matrix inequalities, and that when a solution exists a separation principle holds. That is, a stabilizing output-feedback controller can be found. The procedure is demonstrated using a benchmark problem of nonlinear control: the Moore-Greitzer jet engine compressor model.
\end{abstract}

\section{Introduction}

Output-feedback control design for nonlinear systems remains a challenging problem, because of the unlimited diversity of possible nonlinearities \cite{slotine1991applied},  \cite{isidori1995nonlinear}, \cite{kokotovic2001constructive}. In this paper we present new results on output-feedback design based on control contraction metrics \cite{manchester2014control}, \cite{manchester2014control2} that recover many attractive properties of linear control design, including separation of controller and observer, and convexity of optimization problems.

State-feedback controllers for nonlinear systems can be characterized by the existence of control Lyapunov functions \cite{sontag1983lyapunov}, \cite{artstein1983stabilization}, however these may be difficult to find \cite{rantzer2001dual}. Constructive methods, such as feedback linearization \cite{isidori1995nonlinear}, backstepping \cite{krstic1995nonlinear}, and energy-based methods  \cite{schaft1999l2} are generally applicable only to a limited class of systems. Nonlinear MPC is emerging as a feasible tool (see, e.g. \cite{diehl2009efficient}) but despite some clear benefits, it generally remains difficult to predict or analyse performance of nonlinear MPC schemes by any method other than exhaustive simulations.

Recently there has been signficant interest in using methods of convex optimization to search numerically for feedback control systems. Particular formulations include including density functions \cite{rantzer2001dual}, \cite{prajna2004nonlinear}, occupation measures \cite{lasserre2008nonlinear}, and LQR-Trees \cite{tedrake2010lqr}. These methdos all produce a  state-feedback controller.

There are several approaches to observer design for nonlinear systems, including the circle criterion, \cite{arcak2001nonlinear}, high-gain observers \cite{khalil2002nonlinear}, sliding mode \cite{slotine1991applied}, and contraction theory \cite{Lohmiller98}. However, unlike linear systems, even if a state-feedback controller and a stable observer are found, they can not necessarily be combined to give a stabilizing output-feedback control system. The particular properties of systems, controllers, and observers must be analysed to ensure output feedback stability \cite{atassi2000separation}, \cite{teel1994global}, \cite{arcak2005certainty}.

Contraction analysis \cite{Lohmiller98}, \cite{wang2005partial}, is based on the study of {\em differential} dynamics. Roughly speaking, if {\em all} solutions of a nonlinear system are locally stable, then all solutions converge. Thus global stability results are derived from local criteria, and the problem of motion stability is decoupled from the choice of a particular solution. The search for a contraction metric can be formulated as a convex optimisation problem using sum-of-squares programming  \cite{aylward2008stability} and can be extended to the study of limit cycles \cite{manchester2014transverse}.

A contraction metric can be thought of as a Riemannian metric with the additional property that differential displacements get smaller (with respect to the metric) under the flow of the system. A control contraction metric has the property that differential displacements can be {\em made} to get shorter by control action. This is analogous to the relationship between a Lyapunov function and a control Lyapunov function.

In \cite{manchester2014control} we derived a procedure for state-feedback control design in the form of state-dependent linear matrix inequalities. This is an attractive property because it opens the door to solution via convex optimization methods, such as sum-of-squares \cite{parrilo2000structured}. In \cite{manchester2014control2} we discussed the duality of this control design method with a type of metric observer design, similar to that proposed in \cite{sanfelice2012convergence}. 

In this paper, we show that such dual designs can be used to construct stabilizing output-feedback controllers. That is, a so-called ``separation principle'' holds for controllers and observers based on control contraction metrics. We illustrate this method on a classic benchmark problem in nonlinear control: the Moore-Greitzer jet engine system.


\section{Preliminaries}

For most of this paper, we will consider a nonlinear time-dependant control-affine system
\begin{equation}\label{eq:sys}
\dot x(t) = f(x(t),t) +B(t)u(t), \ y(t) = C(t)x(t)
\end{equation}
where $x(t)\in\RR^n, u(t)\in\RR^m, y(t)\in\RR^p$ are state, control input, and measured output, respectively, at time $t\in\RR^+:=[0,\infty)$. The function $f: \RR^n\times \RR^+ \rightarrow \RR^n$ is assumed to be smooth, and $B(t)$ and $C(t)$ are possibly time-dependent matrices of appropriate dimension.

Contraction analysis is the study of \eqref{eq:sys} by way of the associated system of differential dynamics:
\begin{equation}\label{eq:diffdyn}
\dot\delta_x(t) = A(x,t)\delta_x(t)+B(t)\delta_u(t), \ \delta_y(t) = C(t)\delta_x(t)
\end{equation}
where $A(x,t) = \pder{x}f(x,t)$ is the Jacobian matrix.

In \cite{manchester2014control} we also considered more general systems $\dot x= f(x,u,t)$ but for the moment we note that many systems not naturally appearing in the form \eqref{eq:sys} can be put in that form, either exactly or approximately, by change of variables or introducing new states.

In this paper we study dynamic output-feedback control-systems for \eqref{eq:sys}. Such a controller is a system of the form
\begin{equation}
\dot x_c(t) = f_c(x_c(t), y(t), t), \ u(t) = g_c(x_c(t), y(t), t).
\end{equation}
Frequently the controller will take the form of a state-estimator combined with a state-feedback controller.

The main objective is to design such a feedback system so that the behaviour of the closed-loop system
\begin{align}\label{eq:clsys}
\dot x(t) &= f(x(t),t)+B(t)g_c(x_c(t), C(t)x(t),t) \\
\dot x_c(t) &= f_c(x_c(t), C(t)x(t), t),
\end{align}
is stable.

In particular, given a target trajectory $(x^\star, u^\star)$, feasible for \eqref{eq:sys},  it is desired that there exists positive constants $K$ and $\lambda$ such that \[|x(t)-x^\star(t)|\le K e^{-\lambda t}|x_0-x^\star(0)|\] for all $x(0)$, and the controller state $x_c(t)$ remains bounded (assuming the target trajectory itself remains bounded).

Following \cite{manchester2014control}, we say a system is {\em universally exponentially stabilizable} by output feedback if {\em every} feasible solution $(x^\star, u^\star)$ is globally exponentially stabilizable. This is a stronger condition than global stabilizability of a particular solution (e.g. the origin).

\section{Control Contraction Metrics}\label{sec:ccm}

We now give the basic idea of a control contraction metric (CCM), proposed by the authors in \cite{manchester2014control}. Suppose a system has the property that every solution is locally stabilizable, i.e. the time-varying linear system \eqref{eq:diffdyn}
is stabilizable, where the partial derivatives are evaluated along any particular solution $x(t), u(t)$ of \eqref{eq:sys}. Each local controller may have small region of stability, but if a ``chain'' of states joining the current state $x$ to $x^\star(t)$ is stabilised, in the sense that if each ``link'' in the chain gets shorter, then $x(t)$ is driven towards $x^\star(t)$.

Construction of a CCM is based on taking this concept to the limit as the number of links in the chain goes to infinity, and becomes a smooth path $\gamma(s)$ connecting $x^\star(t)$ and $x(t)$ in the state space. The differential dynamics \eqref{eq:diffdyn} describe the dynamics of infinitesimal path segments. Now, suppose one can find a Riemannian metric $V(x,\delta)=\delta'M(x,t)\delta$ which verifies that a differential feedback law $\delta_u = K(x,t)\delta_x$ is stabilising, i.e.
\begin{align}
(A+BK)'M+M(A+BK)+\dot M<0
\end{align}
then we refer to $M(x)$ as a control contraction metric. The right hand side above can be replaced with $-2\lambda M$ for exponential stability with rate $\lambda$, or various other supply rates for other design criteria such as differential passivity \cite{manchester2014transverse}, \cite{vdsNolcos}, \cite{ForniNolcos} or small gain \cite{manchester2014control2}. The control signal applied is then computed by integrating the differential control signals $\delta_u$ along the path $\gamma$, i.e. 
\begin{equation}\label{eq:control_int}
u(t) = u^\star(t)+\int_\gamma K(\gamma(s))\pder[\gamma]{s}ds.
\end{equation}
Roughly speaking, one stabilises a smooth nonlinear system by stabilising an {\em one-parameter family of linear systems}.

In a sense, this is a generalisation of the concept of a control Lyapunov function to differential dynamics. The advantage is that the rich repertoire of design techniques for linear systems using linear matrix inequalities (LMIs) can be adapted to nonlinear systems as {\em pointwise}-LMIs. Furthermore, pointwise LMIs are now computationally tractable for many important systems thanks to recent advances in semialgebraic optimisation \cite{parrilo2003semidefinite}, \cite{Blekherman13}.

The following definition is central to this paper:
\begin{defin}A function $V(x, \delta_x, t)=\delta_x'M(x,t)\delta_x$, with $\alpha_1I\le M(x,t) \le \alpha_2I$ for some $\alpha_2\ge \alpha_1 >0$, is said to be a {\em control contraction metric} for the system \eqref{eq:sys} if $\pder[V]{x}B(t)=0$ and 
\begin{equation}\label{eq:ctrbmet}
\pder[V]{\delta_x}B(t)=0 \Longrightarrow \pder[V]{t} + \pder[V]{x}f(x,t)+\pder[V]{\delta_x}A(x,t)\delta_x<0
\end{equation}
for all $x, t$.
\end{defin}
Non quadratic Riemann-Finsler metrics can be used without substantial changes to the theory -- see \cite{forni2012differential} for a thorough exploration in the context of stability analysis -- but the associated computational problems for control design are more computationally challenging. In this paper we restrict attention to Riemannian metrics.

We will also make use of a Riemannian distance function between any two points at a given time $d(x_1, x_2, t): \RR^n \times \RR^n \times \RR^+ \rightarrow \RR^+$ defined like so: let $\Gamma(x_1, x_2)$ denote the set of all smooth paths connecting $x_1$ and $x_2$, where each $\gamma\in\Gamma(x_1,x_2)$ is parametrised by $s\in [0, 1]$, i.e. $\gamma(s):[0, 1] \rightarrow \RR^n$. The path length of $\gamma$ is then defined as
\[
L(\gamma,t) := \int_0^1 D\left(\gamma(s),\frac{\partial}{\partial s}\gamma(s),t\right)ds
\]
where $D(x, \delta_x, t) = \sqrt{\delta_x'M(x,t)\delta_x}$ is the length of a differential line element $\delta_x$ with respect to the metric $M(x,t)$. The distance between two points is then defined as
\[
d(x_1, x_2, t) = \min_{\gamma\in\Gamma(x_1, x_2)}L(\gamma,t)
\]
The existence of a minimizing path, which we denote $\gamma_{x_1}^{x_2}(t,s)$,  is implied by the Hopf-Rinow Theorem \cite{spivak1973comprehensive}.


We now briefly summarize some of the results of \cite{manchester2014control} that we will use in this paper.

\begin{theorem}\cite{manchester2014control} If a control contraction metric exists for a system of the form \eqref{eq:sys}, then the system is universally stabilizable by static state feedback.
\end{theorem}

A convex criterion for exponential stabilization is the following:
\begin{theorem}\cite{manchester2014control}
Consider the system \eqref{eq:sys} with differential dynamics \eqref{eq:diffdyn}. If there exists a matrix function $W_c(x,t)\in S_+^n$, $\rho(x,t)\ge 0$ and  $\alpha_2\ge\alpha_1>0$ such that
\begin{align}
\alpha_1I\le W_c\le \alpha_2 I, \label{eq:Wcond} \\
-\dot W_c + W_cA' + AW_c -\rho BB' &\le -2\lambda W_c,\label{eq:Wdotcondexp}
\end{align}
for all $x, u, t$, then the system is universally exponentially stabilizable with rate $\lambda$ by state feedback. In particular, $V(x,\delta,t) = \delta'M_c(x,t)\delta$ is a control contraction metric with $M_c(x,t) = W_c(x,t)^{-1}$, and the following differential gain is stabilizing: $K(x,t) = -\frac{1}{2} \rho(x,t) W(x,t)^{-1}B'$ when used with \eqref{eq:control_int}.
\end{theorem}

In the above condition, $\dot W_c(x,t)$ is a matrix with the $i, j$ element given by $\pder[w_{i,j}]{t}+\pder[w_{i,j}]{x}(f(x,t)+B(t)u)$, where $w_{i,j}$ is the $i,j$ element of $W_c(x,t)$.

\section{Observer Contraction Metric}

It is well-known that the problems of control design and observer design for linear systems have a very attractive ``duality'' (see, e.g., \cite{hespanha2009linear}). It was recently shown by the authors that such a relation also holds for designs based on contraction metrics, building upon past work of \cite{Lohmiller98} and \cite{sanfelice2012convergence}. 

We will call a nonlinear system {\em universally detectable} if the following condition holds: any two solutions $x_1(t), x_2(t)$ that induce identical outputs $y_1(t)=y_2(t) \ \forall t$ have the property that $x_1(t) \rightarrow x_2(t)$ as $t\rightarrow \infty$. I.e. indistinguishable states are convergent.

For such systems, It has been shown by the authors that the existence of a matrix function $W_o(x,t)$ and scalar function $\rho(x,t)\ge 0$ satisfying the  following condition
\begin{equation}\label{eq:OCM}
\dot W_o+A'W_o+W_oA-\rho C'C\le -2\lambda W_o,
\end{equation}
with $W_o(x,t)$ bounded above and below as in \eqref{eq:Wcond}, guarantees existence of an exponentially stable observer \cite{manchester2014control2}.

We note that a similar condition has appeared before in the literature giving necessary conditions for existence of a particular class of observer \cite{sanfelice2012convergence}. By Finsler's theorem, \eqref{eq:OCM} is equivalent to the statement that
\[
C(t)\delta_x(t)= 0 \Longrightarrow \frac{d}{dt}[\delta_x'W_o(x,t)\delta_x]\le -\lambda \delta_x'W_o(x,t)\delta_x.
\]
That is, in directions orthogonal (with respect to $W_o$) to the subspace spanned by the columns of $C(t)$, the system is contracting. Note that in this conditions $W(x,t)$ is the contraction metric, whereas in \eqref{eq:Wdotcondexp} $M(x,t)=W_c(x,t)^{-1}$ was the contraction metric.

\subsection{Construction of an Observer}\label{sec:obs}
Suppose condition \eqref{eq:OCM} is satisfied, and choose an initial state estimate $\hat x(0)$. At each time $t$, define the set $\mathfrak X_y(t):=\{x:C(t)x=y(t)\}$ as the set of states perfectly consistent with the measurement at time $t$. Now let $\gamma(s,t)$ be the shortest path, with respect to the Riemannian metric $\delta'W_o(x,t)\delta$, between $\hat x(t)$ and the set $\mathfrak X_y(t)$.

Then construct an observer with the following dynamics:
\begin{equation}\label{eq:Obs}
\dot {\hat x}(t) = f(\hat x(t),t)+\int_0^T  K(\gamma(s))\delta_y(s)ds
\end{equation}
where
\[
K(x)=\frac{1}{2}\rho(x,t)W(x,t)^{-1}C(t)'.
\]
Note that if $W$ and $\rho$ are independent of $x$, this reduces to a standard Luenberger-type observer
\[
\dot {\hat x}(t) = f(\hat x(t),t)+ K(t)(y(t)-C(t)\hat x(t)).
\]

A set $S$ is called geodesically convex with respect to a metric $\delta'W(x,t)\delta$ if for any two points $x_1$ and $x_2$ in $S$, the minimal geodesic connecting them remains in $S$. The following result was recently proved by the authors:
\begin{theorem} \cite{manchester2014control2} If the condition \eqref{eq:OCM} hold and, additionally, the set $\mathfrak X_y(t)$ is geodesically convex with respect to the metric $W(x,t)$, then the observer constructed above converges, i.e $\hat x(t) \rightarrow x(t)$. By construction, the system is universally detectable.
\end{theorem}

Note that if $W$ is independent of $x$, then any set of the form $\{x:y(t) = C(t)x\}$ -- i.e. an affine variety -- is geodesically convex.

\section{Separation of Observer and Controller Design}

If a linear system is both stabilizable and detectable, then it is stabilizable by output feedback using a combination of a stable observer and a state-feedback controller. This extremely useful property, known as the ``separation principle'', fails to hold for general nonlinear systems.

In this section we prove the main theoretical result of this paper: that if the stronger property holds that a system is {\em universally stabilizable and detectable}, then a separation principle does indeed hold.

Firstly, we prove a result on input-to-state stability of systems controlled using CCMs:

\ 

\begin{theorem}\label{thm:expbound}
Consider a system of the form
\[
\dot x = f(x)+Bu+w
\]
where $x$ and $u$ are state and control input, and $w$ is a disturbance input. If the state-feedback control design proposed in Section \ref{sec:ccm} is used then the following holds:
\[
\dot d(t) \le -\lambda d(t) + \sqrt{\alpha_1} \|w\|
\]
where $d(t)$ is the Riemannian distance from $x(t)$ to $x^\star(t)$, with respect to the metric $M(x)$ and $\Theta(x)'\Theta(x) = M(x)$.
\end{theorem}

\ 

\begin{proof}
The distance at time $t$ is
\[
d(x(t), x^\star(t)) = \int_\gamma \|\Theta(x)\delta_x\|ds
\]
with $x = \gamma(s)$ and $\delta_x = \pder[\gamma]{s}$.

The differential dynamics at each point on the geodesic $\gamma(s)$, $s\in[0,1]$ satisfy
\[
\dot\delta_x = (A+BK)\delta_x + w.
\]
where $\delta_x = \pder[\gamma]{s}$. At each point on the geodesic, we consider the differential change of coordinates $\delta_z = \Theta(x)\delta_x$, and we have
\[
\dot \delta_z = F\delta_z+\Theta(x)w,
\]
where $F$ is a generalised Jacobian \cite{Lohmiller98} satisfying $\delta_z'F\delta_z\le -\lambda \|\delta_z\|$ for all $\delta_z$.

The control contraction condition implies that
\[
\delta'(\dot M + (A+BK)'M + M(A+BK) + 2\lambda M)\delta \le 0,
\]
which gives the following expression for the derivative of the length of a differential line element $\sqrt{\delta'M\delta} = \sqrt{\delta_z'\delta_z}=\|\delta_z\|$:
\[
\frac{d}{dt} \sqrt{\delta_z'\delta_z} = \frac{\delta_z'F\delta_z + \delta_z'\Theta(x)w}{\sqrt{\delta_z'\delta_z}}.
\]
So by the bound on the generalized jacobian, and by applying the Cauchy-Schwarz inequality to give $\delta_z'\Theta(x)w\le \|\Theta(x)w\|\|\delta_z\|$, we have
\[
\frac{d}{dt} \|\delta_z\| \le -\lambda \|\delta_z\|+\|\Theta(x)w\|
\]
From the bound $W_c\ge \alpha_1 I$ we have $M_c\le \alpha_1I$, and since $M = \Theta'\Theta$ we have$\|\Theta(x)w\|\le \sqrt{\alpha_2}\|w\|$ for all $x$ on the path $\gamma_t$, and so integrating the above inequality along $\gamma_t$ gives the result of the theorem.
\end{proof}

We are now ready to state the main theoretical result of the paper.

\ 

\begin{theorem}
A system of the form \eqref{eq:sys}. 
Construct the observer as in Section \ref{sec:obs} giving the state estimate $\hat x(t)$,  and construct the control signal as in Section \ref{sec:ccm}, but with $\hat x(t)$ in place of $x(t)$, then the closed-loop system is exponentially stable.
\end{theorem}

\

\begin{proof}
The closed-loop dynamics can be written as
\[
\dot x = f(x,t) + Bk(x,t) + B(k(\hat x,t) - k(x,t))
\]
where $k(x,t) = u^\star(t)+\int_\gamma K(x,t)\delta_xds$ is the control law from Section \ref{sec:ccm}.

The the exponential convergence of the state estimates $\hat x\rightarrow x$ and smoothness of $k$ implies that $k(\hat x,t) - k(x,t)$ is bounded and converges to zero asymptotically. This further ensures boundedness of the state $x$, by Theorem \ref{thm:expbound} and the boundedness of $\hat x^\star$. On any compact set $K(x,t)$ is uniformly bounded, and so we can also affirm that $k(\hat x)-k(x)\rightarrow 0$ exponentially. 

The uniform boundedness of $\Theta$ then guarantees the existence of some $\alpha>0, \beta>0$ such that
\[
\|\Theta(\gamma_t(s)) B(k(\hat x(t)) - k(x(t)))\| \le \beta e^{-\alpha t}
\]
for all $s\in [0,1]$ and $t\ge 0$.

Now, using again Theorem \ref{thm:expbound} with $w(t) = B(k(\hat x,t) - k(x,t))$ it follows that 
\begin{align}
d(t) \le& d(0)e^{-\lambda t} +\notag \\
& \int_0^T e^{-\lambda(T-t)}\  \max_{s\in [0,1]} \|\Theta(\gamma_\tau(s)) B(k(\hat x,t) - k(x,t))\| dt\notag \\
&\le d(0)e^{-\lambda t} +  \beta\int_0^T e^{-\lambda(T-t)} e^{-\alpha t}dt\notag
\end{align}
which implies that $d(t) \rightarrow 0$ exponentially. 

The uniform boundedness of the metric then implies that $\|x(t) - x^\star (t)\|$ converges to zero exponentially.
\end{proof}

\section{Simplification for Constant Metrics}

We have presented a general construction of output-feedback controller that allows state-varying contraction metrics for both the observer and controller parts. In the special case that one can find constant metrics -- i.e. independent of $x$ this construction can be substantially simplified.

For a constant metric $M_c=W_c^{-1}$, the geodesic joining $\hat x$ and $x^\star$ is always a straight line. Defining $\Delta_c(t) = x^\star(t)-\hat x(t)$ as the error between the current state estimate and the target trajectory, the control signal \eqref{eq:control_int} is:
\[
u(t) = u^\star(t)+\frac{1}{2}\int_0^1 \rho(x(s,t))M_cB' \Delta_c(t) ds
\]
where $x(s,t) = \hat x(t) + s\Delta_c(t)$. But since $M_c, B$ and $\Delta_c$ are independent of $s$, they can be taken outside the integral and we get
\[
u(t) = u^\star(t)+\frac{1}{2}\left(\int_0^1 \rho(x(s,t))ds \right)M_cB' \Delta_c(t).
\]
Furthermore, the integral in the brackets is just a one dimensional integral of $\rho$, with a particular one-dimensional locus of points substituted as argument. In the case when $\rho$ is a polynomial (as when sum-of-squares programming is used to compute the controller) this can be expressed analytically as a difference of two polynomials in $x$. Hence the state-feedback portion of the control law becomes expressed as a polynomial in $\hat x(t)$.

Similarly, for the observer design, the minimal geodesic joining $\hat x(t)$ to $\mathfrak X_y(t)$ is a straight line and can be analytically constructed, since it is simply a linearly constrained weighted least-squares problem:
\begin{align}
&\bar x(t) = \arg\min_{x\in\mathbb R^n} \ (x-\hat x(t))'W_o(x-\hat x(t) \notag \\
&\textrm{ subject to } y(t) = C(t)x.\notag
\end{align}
This can be obtained analytically be solving the linear system:
\[
\begin{bmatrix} W_o & C(t)' \\ C(t)& 0\end{bmatrix}\begin{bmatrix} \bar x(t) \\ \lambda\end{bmatrix} = \begin{bmatrix} W_o \hat x(t) \\ y(t)\end{bmatrix}.
\]
Then the observer can be constructed in explicit form as for the controller, except with $\Delta_o = \hat x(t) -\bar x(t)$.

Hence, in the case of polynomial dynamics and constant control and observer contraction metrics, and output-feedback control can be explicitly constructed as a polynomial control law and a polynomial state observer.
\section{Application Example: Jet Engine Oscillations}

The Moore-Greitzer model, a simplified model of surge-stall dynamics of a jet engine  \cite{moore1986theory}, has motivated substantial development in nonlinear control design (see, e.g., \cite{moore1986theory},  \cite{krstic1995nonlinear}, \cite{krstic1998useful},  and references therein). In \cite{aylward2008stability}, sum-of-squares programming was applied for robustness analysis of stable solutions, and in \cite{manchester2014transverse} transverse contraction was used to analyse the orbital stability of compressor oscillations. Output feedback control has previously been addressed in \cite{chaturvedi2006output}, \cite{shiriaev2010global}.

A model of surge-stall dynamics was derived in \cite{moore1986theory} based on a Galerkin projection of the PDE on to a Fourier basis. The following reduced model of the surge dynamics is frequently studied:
\[
\begin{bmatrix}\dot\phi \\ \dot\psi \end{bmatrix}  = \begin{bmatrix} -\psi -\frac{3}{2}\phi^2 - \frac{1}{2}\phi^3\\ \phi + u \end{bmatrix}.
\]
with $u$ as the input and a sensor on $\psi$ only. Here $\phi$ is a measure of mass flow through the compressor, and $\psi$ is a measure of the pressure rise in the compressor, under a change of coordinates, see \cite[p. 68]{krstic1995nonlinear}. The source of difficulty is the nonlinearity $-\frac{3}{2}\phi^2 - \frac{1}{2}\phi^3$ which does not satisfy any global Lipschitz bound, and affects the dynamics of the variable $\phi$, which is not directly controlled or measured.

The system can be written in the form \eqref{eq:sys} with $x = [\phi, \ \psi]'$,
\[
f(x) = \begin{bmatrix} -\psi -\frac{3}{2}\phi^2 - \frac{1}{2}\phi^3\\ \phi\end{bmatrix}, B = \begin{bmatrix}0\\ 1\end{bmatrix},C = \begin{bmatrix}0 & 1\end{bmatrix}.
\]
and 
\[
A(x) = \begin{bmatrix} 
-3\phi -\frac{3}{2}\phi^2 & -1 \\ 1 & 0
\end{bmatrix}.
\]

\begin{figure}
\begin{center}
\includegraphics[width=1\columnwidth]{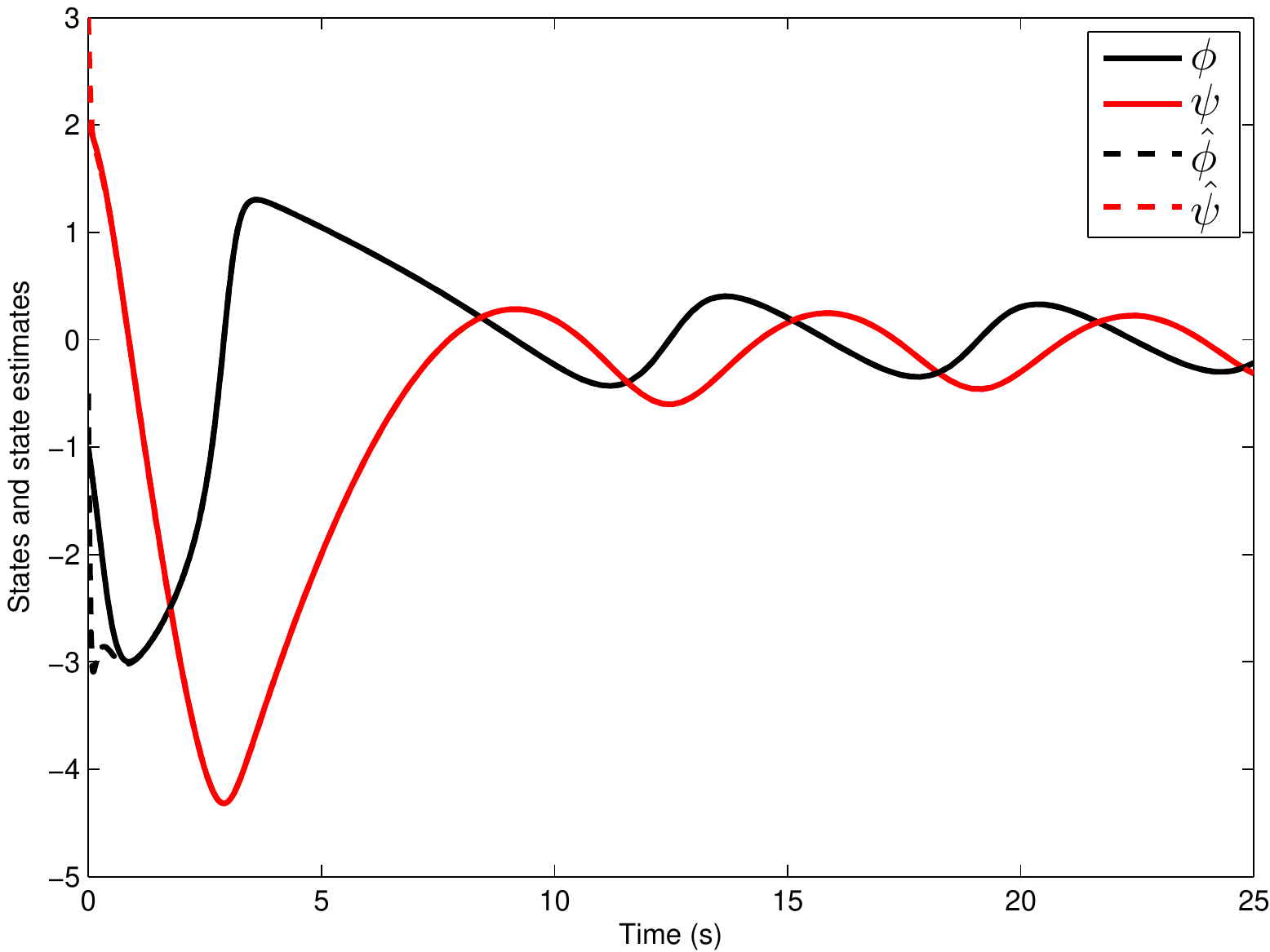}
\caption{States and state estimates of Moore-Greitzer jet engine model in open-loop, exhibiting combustion oscillations.}
\label{fig:OL}
\end{center}
\end{figure}

As is well known, from certain initial conditions this system can exhibit orbitally stable oscillating solutions, as seen in Figure \ref{fig:OL}. This is an inherently nonlinear phenomenon. The objective of the control design is to stabilize these oscillations to zero.

We used Yalmip sum-of-squares module \cite{YALMIP}, \cite{Lofberg09} to set up a two decoupled convex feasibility problems, one for the state-feedback controller:
\begin{align}
-\delta'(W_cA(x)'+A(x)W_c-\rho_c(x)BB' +\lambda W_c) \delta&\in \Sigma^2[x,\delta] \notag\\
\rho_c(x) & \in \Sigma^2[x] \notag\\
W_c\ge \alpha_1I, W_c &\le \alpha_2I, \notag
\end{align}
and one for the observer:
\begin{align}
-\delta' (W_oA(x)'+A(x)W_o-\rho_o(x)C'C + \lambda W_o)\delta& \in \Sigma^2[x,\delta] \notag\\
\rho_o(x) & \in \Sigma^2[x] \notag, \\
W_o\ge \alpha_1I, \ W_o &\le \alpha_2I. \notag
\end{align}
Here the decision variables are the symmetric matrices $W_c, W_o$ and the coefficients of the degree-two polynomials $\rho_c, \rho_o$ -- these appear linearly in the above constraints. The ``abstract'' variables over which the sum-of-squares constraints hold are $x$ and $\delta$, where $\Sigma^2[\cdot]$ denotes a sum of squares constraint with respect to the abstract variable in the argument. 

The bounds $\alpha_1, \alpha_2, \lambda$ were fixed in advance and the same for controller and observer, although they could be different and considered as variables for optimization.

The controller and observer problems were  solved sepearately. The sum-of-squares relaxations each had seven scalar variables, four $2\times 2$ matrix variables, and 27 constraints. Using the commercial solver Mosek version 7 \cite{andersen2000mosek}, each of these semidefinite programs took less than 0.4 seconds to solve on a standard desktop workstation. 

\begin{figure}
\begin{center}
\includegraphics[width=1\columnwidth]{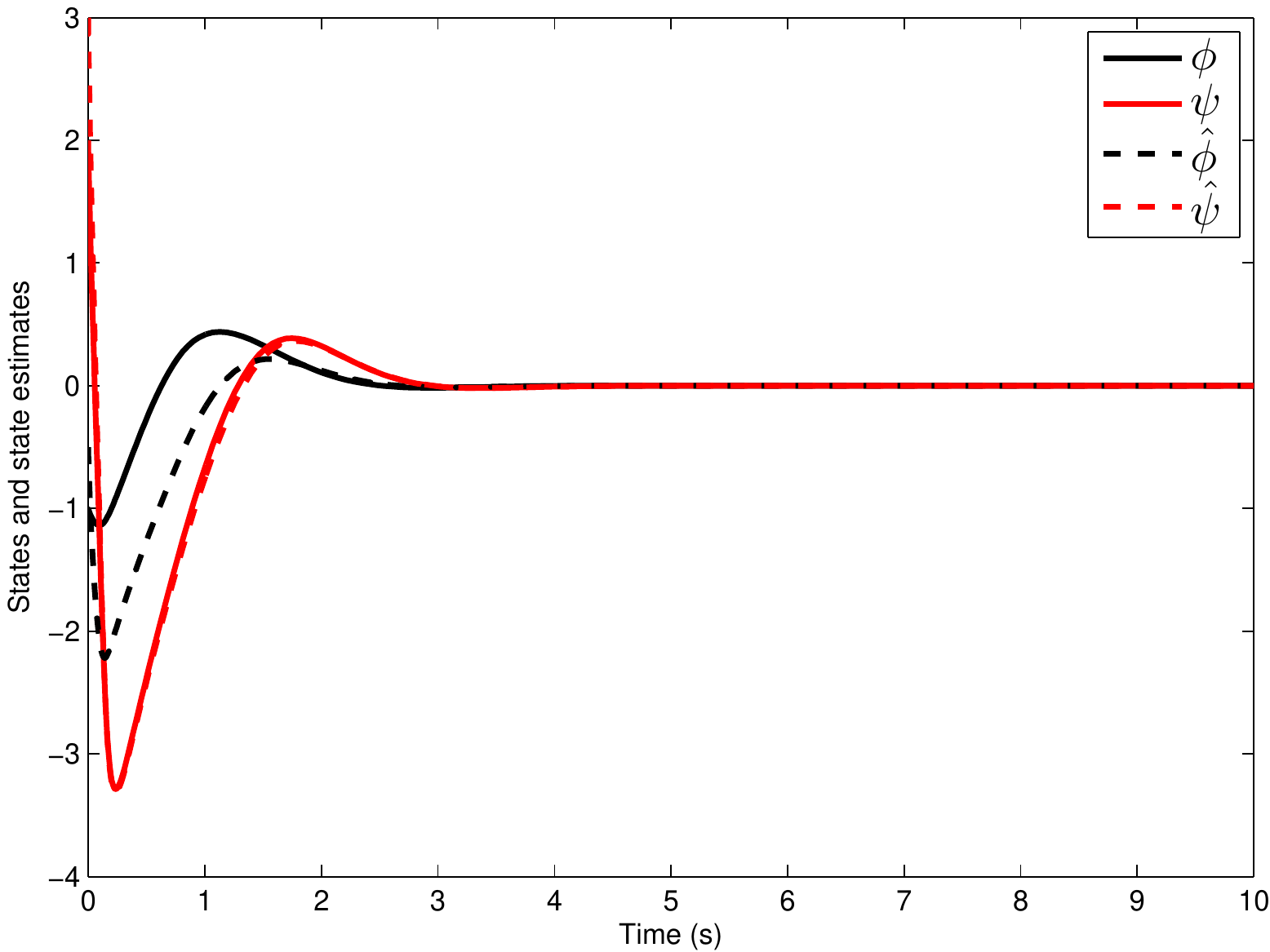}
\caption{States and state estimates of  stabilized  Moore-Greitzer system with convergence rate $\lambda = 0.1, \alpha_1 = 0.1, \alpha_2 = 1.3$.}
\label{fig:SlowResponse}
\end{center}
\end{figure}

The results exhibit the trade-off between speed of convergence and shaping transient response. Specifying a relatively slow rate of convergence of $\lambda = 0.1$, we were able to find a constant metrics with satisfying \eqref{eq:Wcond} with $\alpha_1=0.1$ and $\alpha_2=1.3$. The simulation is shown in Figure \ref{fig:SlowResponse}. As can be seen, the system actually converges significantly faster than this constraint requires. The relative sizes of the upper and lower bounds on $W_c$ and $W_o$ correspond loosely to overshoot: each differential line element in each geodesic (for controller and observer) has strictly forward-invariant sets of the form $\delta'W\delta = c$ for any $c>0$. Hence a $W$ which is ``close'' to identity results in little overshoot, in the system's original coordinates. 

\begin{figure}
\begin{center}
\includegraphics[width=1\columnwidth]{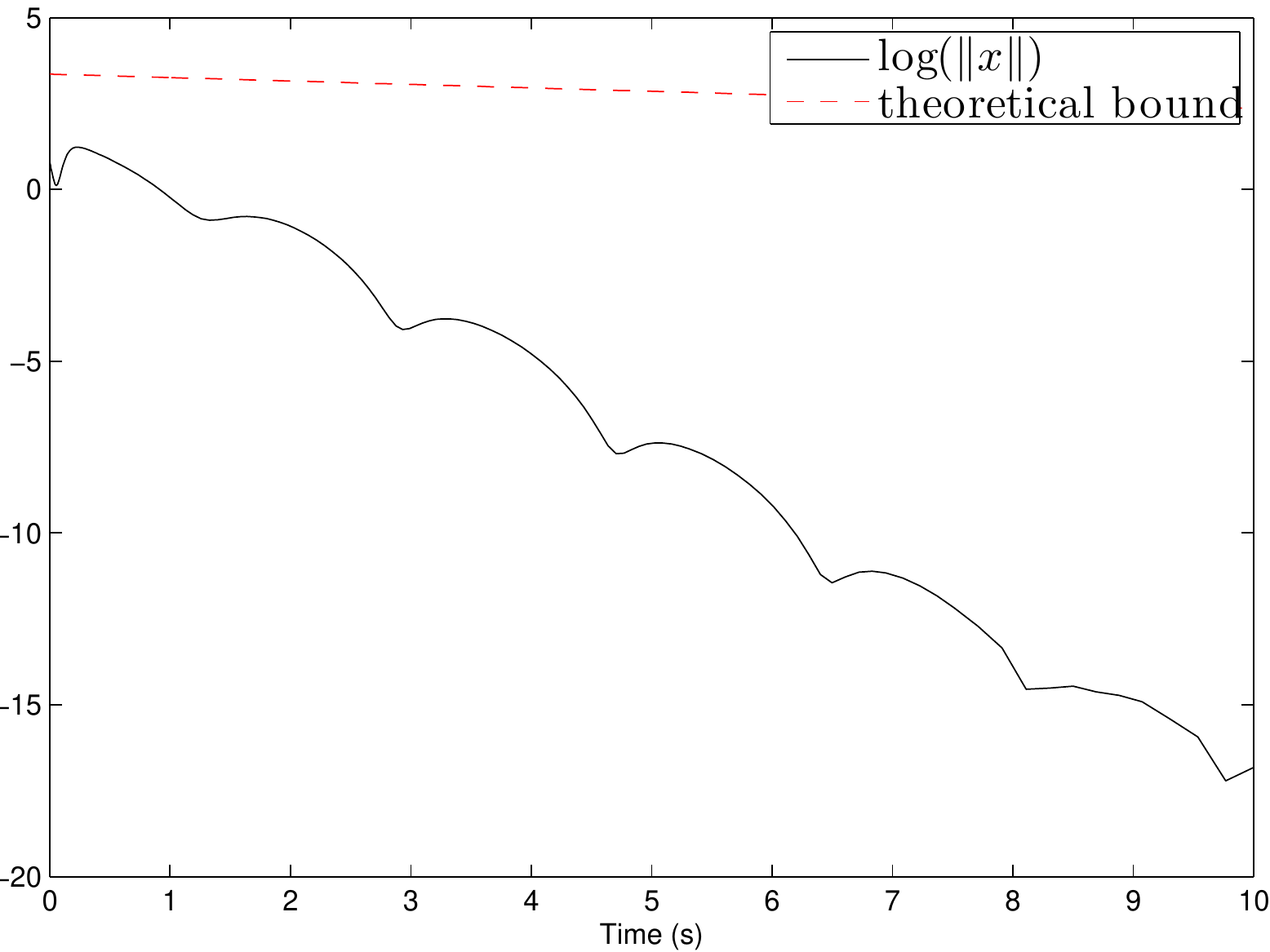}
\caption{Exponential convergence of simulation, and theoretical upper-bound with $\lambda = 1, \alpha_1 = 0.1, \alpha_2 = 1.3$.}
\label{fig:SlowLogBounds}
\end{center}
\end{figure}

\begin{figure}
\begin{center}
\includegraphics[width=1\columnwidth]{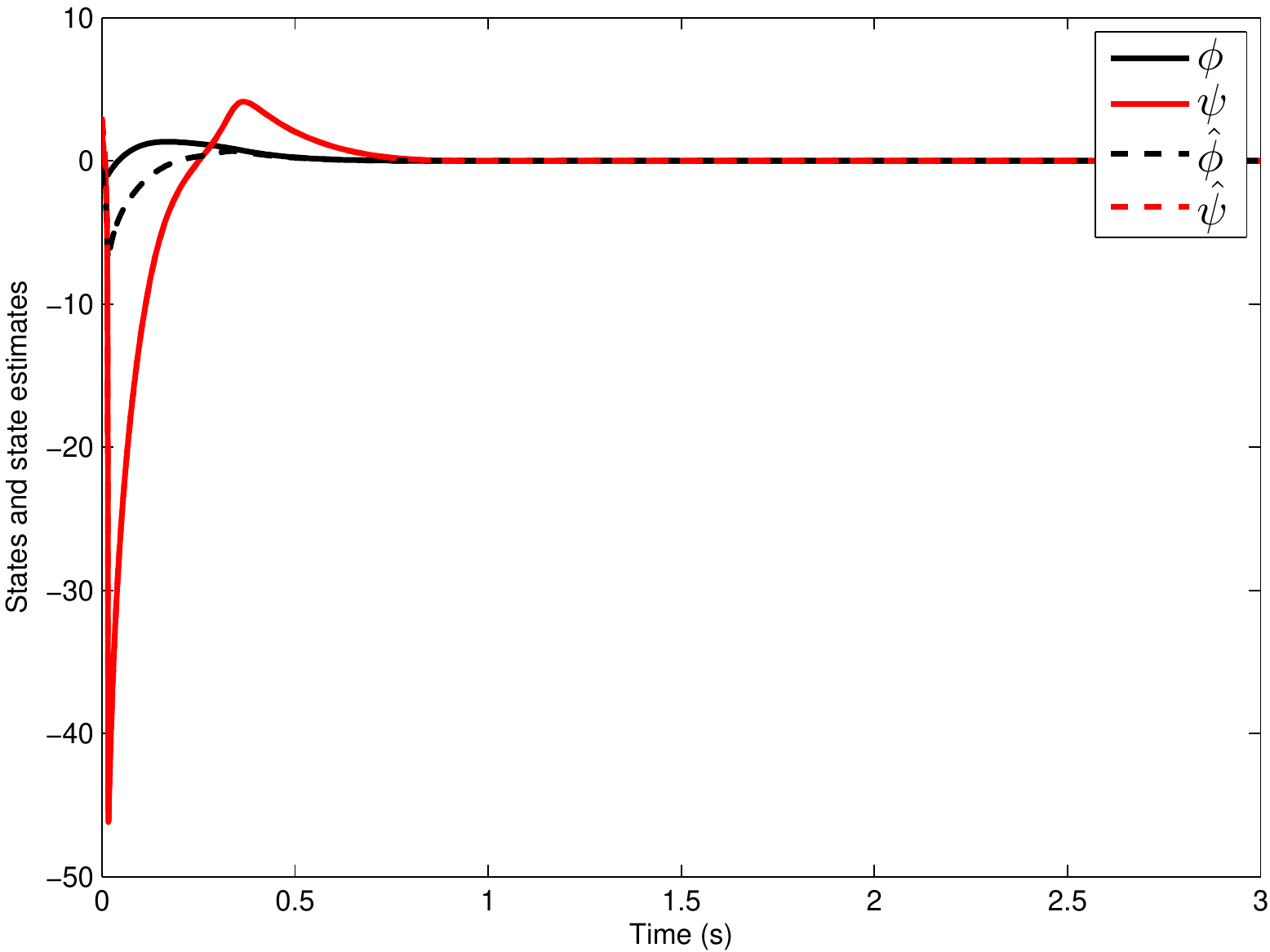}
\caption{States and state estimates of stabilized Moore-Greitzer system with convergence rate $\lambda = 5, \alpha_1 = 0.1, \alpha_2 = 30$.}
\label{fig:MedResponse}
\end{center}
\end{figure}

With a faster rate of convergence, $\lambda = 5$, we needed to expand the range of the bounds on $W_c$  and $W_o$ to $\alpha_1 = 0.1, \alpha_2 = 30$. As can be seen in Figure \ref{fig:MedResponse}, there is some overshoot, or ``peaking'' observed. Note that the scale of the time axis has changed. Such peaking is a common issue in certain high-gain observer designs.

\begin{figure}
\begin{center}
\includegraphics[width=1\columnwidth]{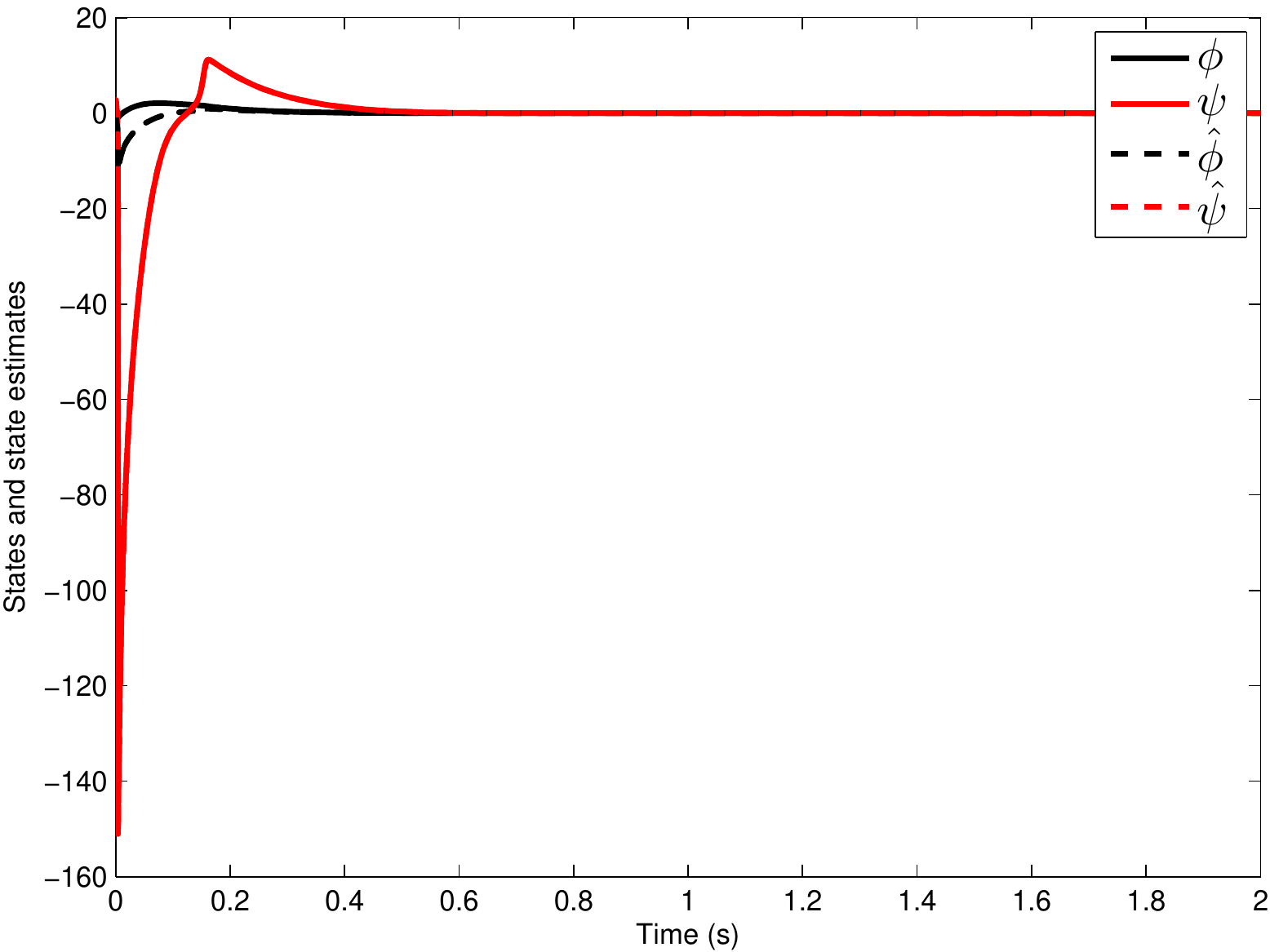}
\caption{States and state estimates of stabilized  Moore-Greitzer system with convergence rate $\lambda = 10, \alpha_1 = 0.1, \alpha_2 = 100$.}
\label{fig:FastResponse}
\end{center}
\end{figure}

\begin{figure}
\begin{center}
\includegraphics[width=1\columnwidth]{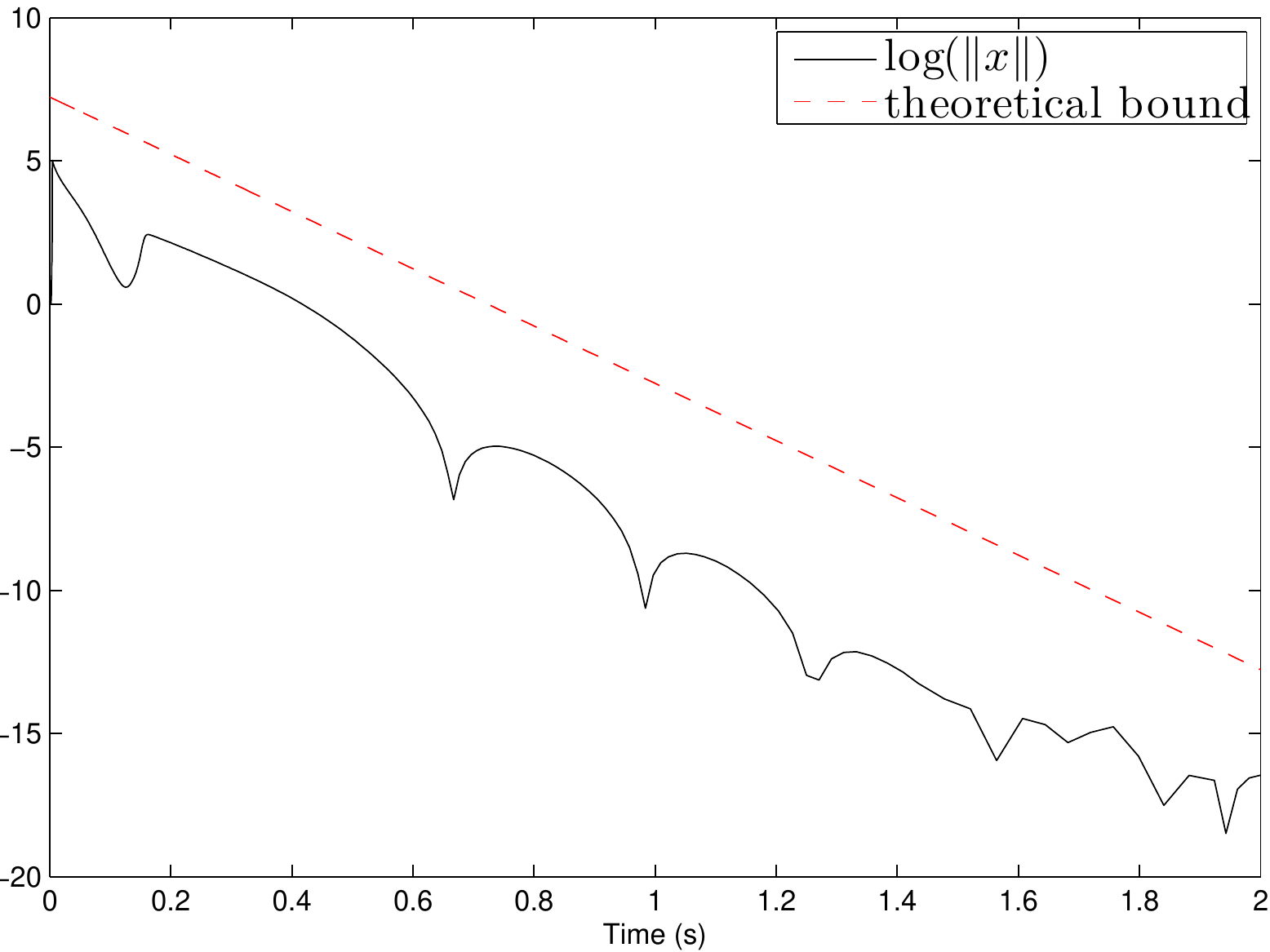}
\caption{Exponential convergence of simulation, and theoretical upper-bound with $\lambda = 10, \alpha_1 = 0.1, \alpha_2 = 100$}
\label{fig:FastLogBounds}
\end{center}
\end{figure}

With an even fast rate of convergence $\lambda = 10$, more severe peaking is observed, see  Figure \ref{fig:FastResponse}. This may be impractical for the real system, so we can see that adjusting the bounds $\alpha_1, \alpha_2$ and $\lambda$ allows a trade-off between convergence and overshoot. Note that the system's convergence rate is quite well approximated by the theoretical bound, as seen by the slopes of the lines in Figure \ref{fig:FastLogBounds}, indicating that in this case the constraints were quite tight.

On the other hand, this is just the behaviour of one particular method of designing a differential feedback gain. An advantage of the control contraction metric formalism is that many techniques from linear control can be directly adapted, e.g. LQG-like output feedback. There is much work be to be done to explore the possibilities.

\begin{figure}
\begin{center}
\includegraphics[width=1\columnwidth]{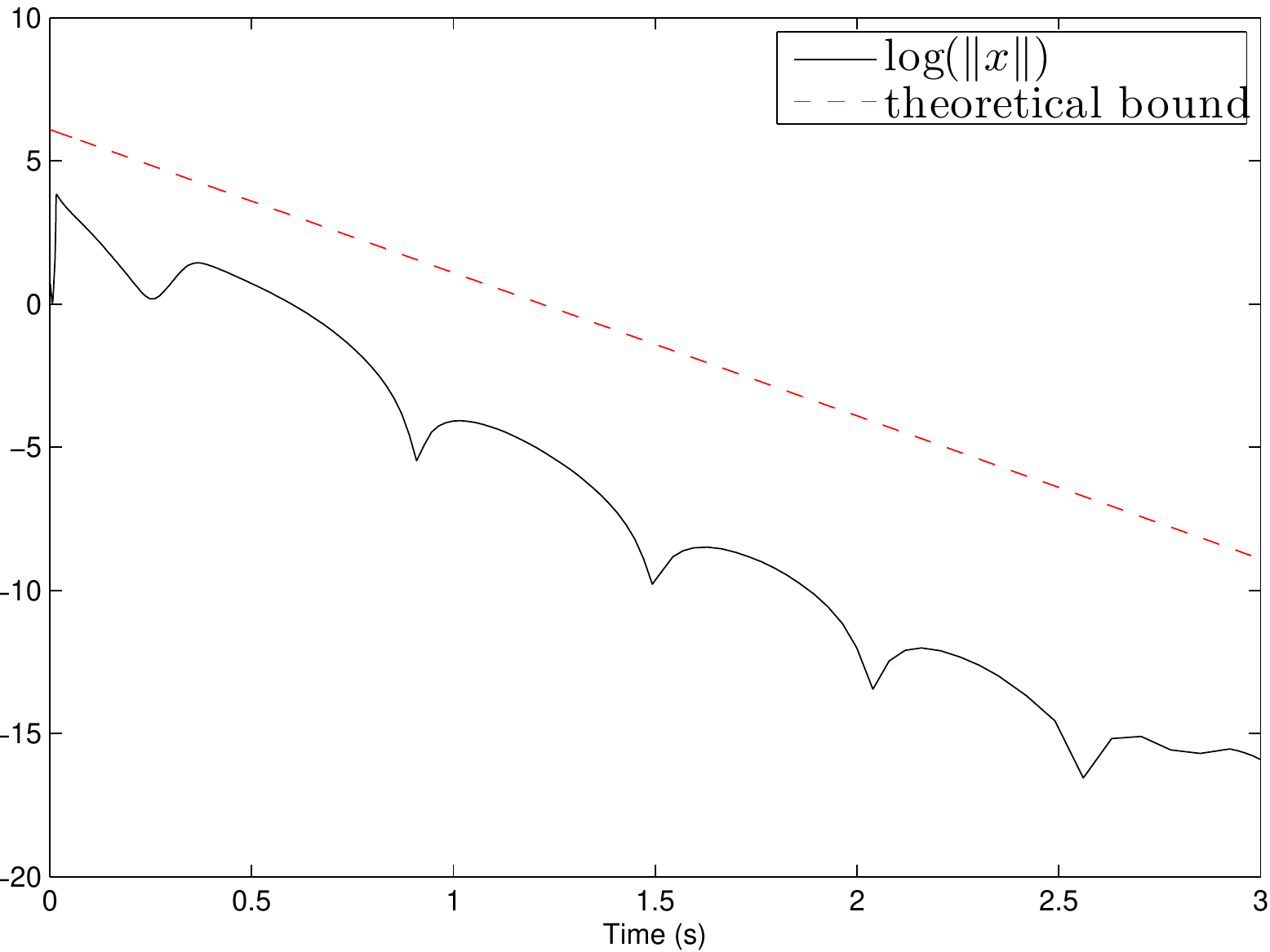}
\caption{Exponential convergence of simulation, and theoretical upper-bound with $\lambda = 5, \alpha_1 = 0.1, \alpha_2 = 30$.}
\label{fig:MedLogBounds}
\end{center}
\end{figure}

\begin{figure}
\begin{center}
\includegraphics[width=1\columnwidth]{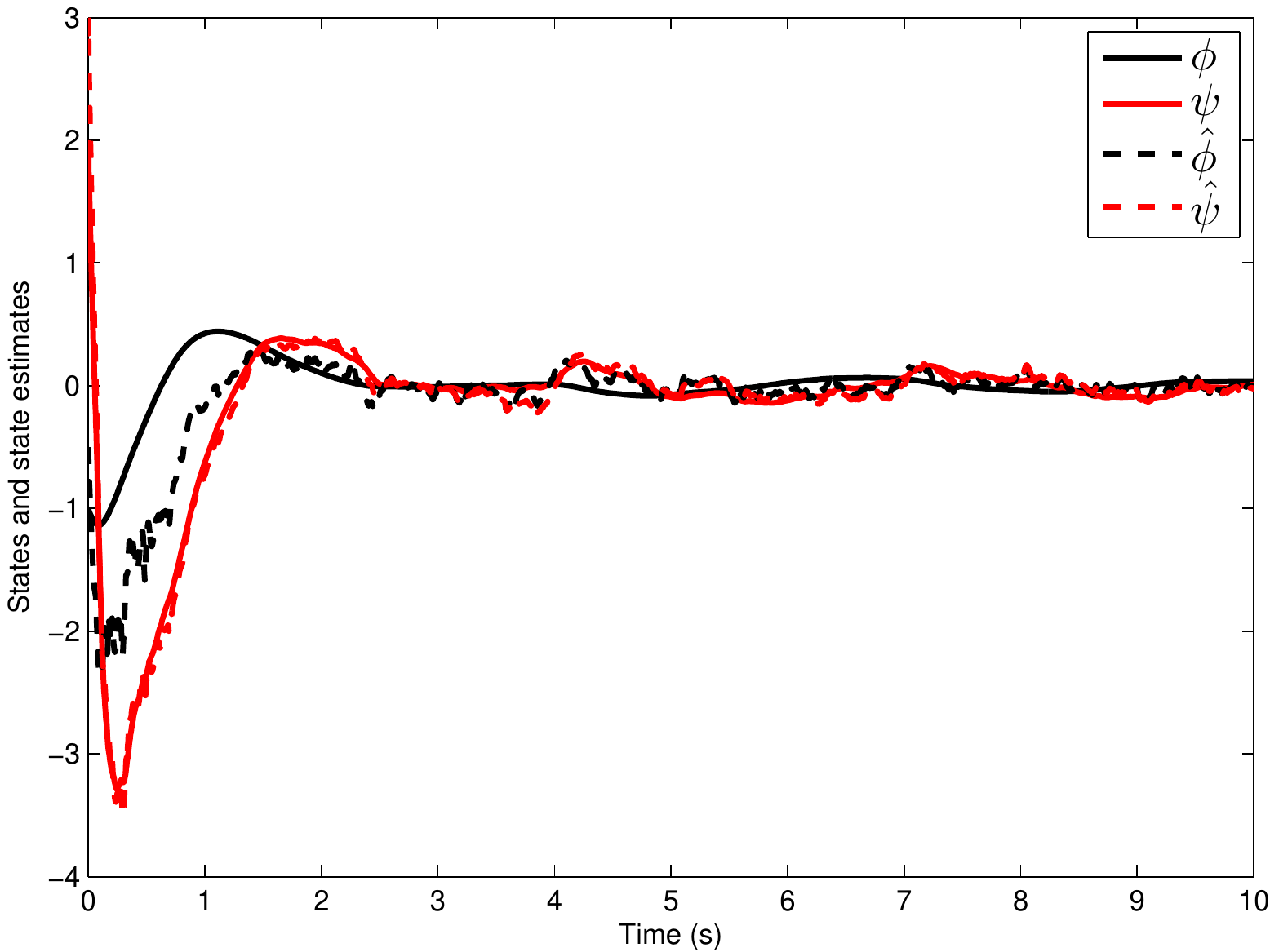}
\caption{States and state estimates of stabilized Moore-Greitzer system with measurement noise.}
\label{fig:NoisyResponse}
\end{center}
\end{figure}

\begin{figure}
\begin{center}
\includegraphics[width=1\columnwidth]{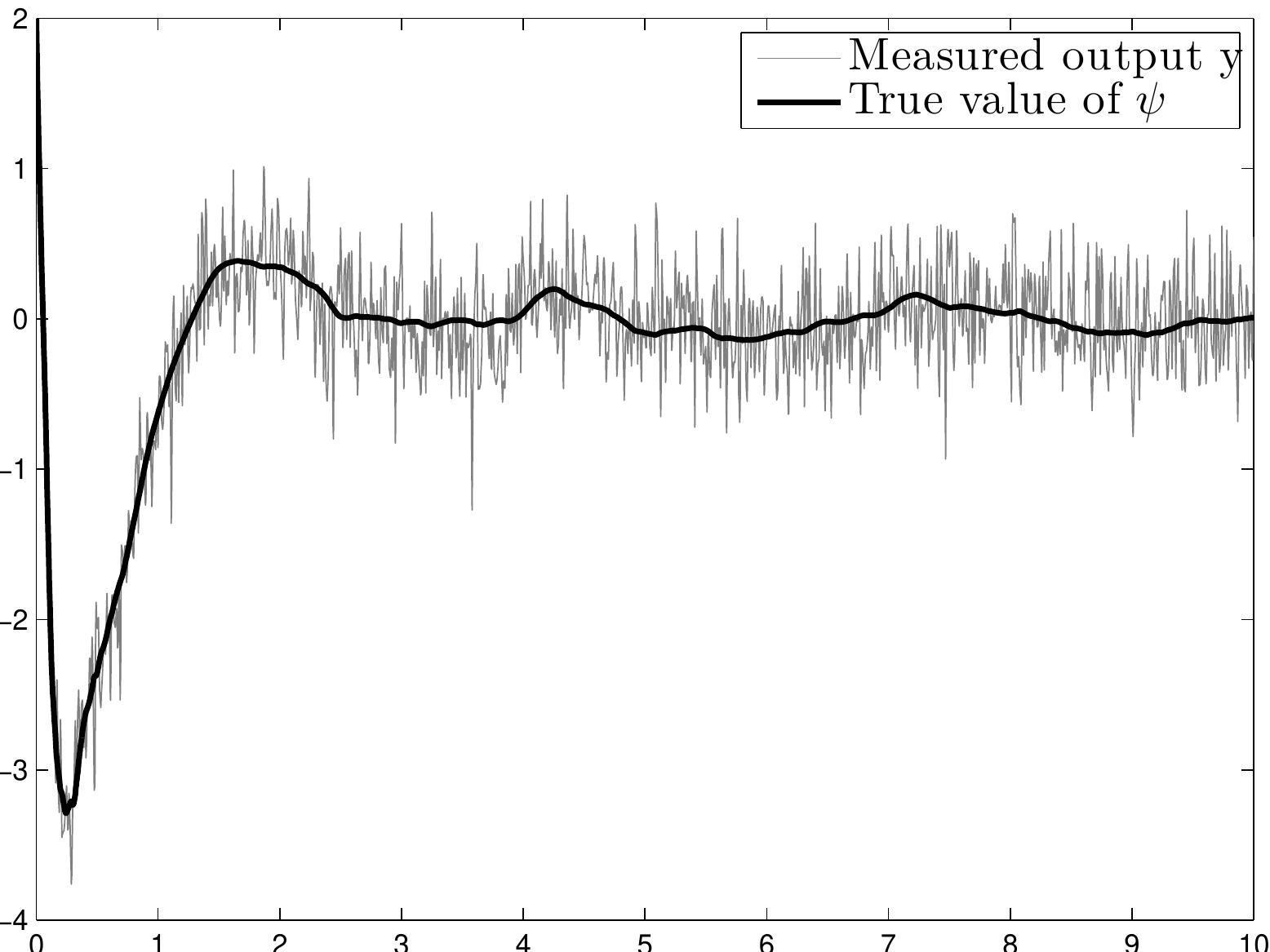}
\caption{Visualisation of measurement noise: measurement $y(t)$ and the true value of $\psi(t)$.}
\label{fig:MeasNoise}
\end{center}
\end{figure}

An output-feedback controller should also be able to operate with noisy measurements. We simulated the closed-loop system with $\lambda = 0.1$ with a measurement noise on the output $y(t) = \psi(t) + n(t)$ where $n(t)$ was a Gaussian white-noise process with standard deviation 0.3. As can be seen from Figure \ref{fig:NoisyResponse}, the controller is very able to keep the system stabilised, despite the quite small signal to noise ratio from the sensor -- see Figure \ref{fig:MeasNoise}.

\section{Conclusions}

In this paper we have shown that recent methods for state-feedback controller design and observer design based on control contraction metrics can be combined to form a stabilizing output-feedback controller. This extends the well-known separation theorem for linear systems to a class of nonlinear systems.

The constructions are computationally tractable, taking the form of pointwise linear matrix inequalities. The numerical search for such metrics forms a convex optimization (or feasibility) problem. When the system has polynomial dynamics, this can be solved using sum-of-squares programming and a semidefinite programming solver, of which there are several free and commercial options.

In future work we will continue to investigate practical examples, especially in robotics, as well as different methods of generating gains. We will also investigate the potential to extend these results to robust and adaptive control design, using the recently developed concept of differential passivity \cite{manchester2014transverse}, \cite{vdsNolcos}, \cite{ForniNolcos}, and differential $L^2$ gain \cite{manchester2014control2}. We will also explore the use of nonlinear convex parameterizations of stable models proposed in \cite{tobenkin2010convex}, \cite{manchester2011}.

\bibliographystyle{IEEEtran}
\bibliography{elib}

\end{document}